\documentclass[a4paper,11pt]{article}
\usepackage[english]{babel}
\usepackage{graphicx}
\usepackage{bbding}
\usepackage[latin1]{inputenc}
\usepackage{fancyhdr}
\usepackage{amsmath,amssymb,color,epsfig}
\usepackage{mathrsfs}
\usepackage{eufrak}
\usepackage{dsfont}
\usepackage{upgreek}
\usepackage{fancyhdr}
\usepackage[dvipsnames]{xcolor}
\usepackage{lipsum}
\newtheorem{theorem}{Theorem}[section]
\newtheorem{proposition}[theorem]{Proposition}
\newtheorem{corollary}[theorem]{Corollary}
\newtheorem{remark}[theorem]{Remark}
\newtheorem{lemma}[theorem]{Lemma}

\def\1{\mathds{1}}

\def\Re{{\rm Re}}

\usepackage{authblk}
\usepackage{hyperref}
\newcommand*{\email}[1]{%
    \normalsize\href{mailto:#1}{#1}\par
    }
\title{A MEAN VALUE INEQUALITIES FOR THE POLYGAMMA AND ZETA
FUNCTIONS}
\author{Mohamed Bouali}
\affil{Department of Mathematics, University of Tunis,\\ University of Tunis El-Manar
  \\ \email{bouali25@laposte.net}}

\begin{document}

\maketitle

\begin{abstract} A recently published result states inequalities of the harmonic mean of the digamma function. In this work, we prove among others results that
 for all positive real numbers $x\neq 1$,
 $$-\gamma<-\gamma H(x,1/x)<\frac{\gamma^2}{\psi\big(H(x,1/x)\big)}<\psi\Big(1/H(x,1/x)\Big)<H\Big(\psi(x),\psi(1/x)\Big),$$
$$H\Big(\zeta(x),\zeta(1/x)\Big)<-2,$$
and for all $x\in(0,1)$
$$\zeta(1/2)<H\Big(\zeta(x),\zeta(1-x)\Big)<-1,$$
$$\frac{\log 4}{1+\log 4}<H\Big(\eta(x),\eta(1-x)\Big)<(1-\sqrt 2)\zeta(1/2).$$
Here, $\psi=\Gamma'/\Gamma$ denotes the digamma function, $\gamma$ is Euler's constant, $\zeta$ is the Riemann's zeta function and $\eta$ is the Dirichlet's eta function.

\end{abstract}

%\section{Introduction}
%\label{sec:introduction}
\section{Digamma function}
 The classical harmonic mean of two real numbers $a$ and $b$ (which are not both equal to
$0$) is defined by
 $$H(a,b)=\frac{2ab}{a+b}.$$
 It satisfies $H(a,1/a)\in(0,1)$ for all $a>0$.

 In a recent paper, Jameson and Alzer \cite{alz3} proved a remarkable mean value inequality for the
digamma function.
$$\psi(x)=\frac{\Gamma'(x)}{\Gamma(x)},$$
where $\Gamma$ is the Euler gamma function.
They showed that for all positive real numbers $x$ the harmonic mean of $\psi(x)$
and $\psi(1/x)$ is greater than or equal to $-\gamma$, that is,
$$-\gamma<H(\psi(x),\psi(1/x)).$$
Equality holds if and only if $x=1$.

A refinement of this inequality is proved in \cite{alz4}. He showed that for $x>0$, $x\neq 1$
$$-\gamma H(x,1/x)<H(\psi(x),\psi(1/x)).$$
The digamma function has interesting applications in numerous fields, like, for instance, the theory of special functions, statistics, mathematical physics and number theory. See \cite{3}, \cite{4}.

The following series and integral representations are valid for $x >0:$
$$\psi(x)=-\gamma-\frac1x+\sum_{k=1}^\infty\frac x{k(x+k)}=\int_0^\infty\Big(\frac{e^{-t}}t-\frac{e^{-x t}}{1-e^{-t}}\Big)dt.$$
$\psi$ has a unique zero $x_0\simeq 1.4616...$  (In what follows, we keep this notation.) Moreover,
$\psi$ is strictly increasing and strictly concave on $(0,+\infty)$ and satisfies the limit relations
$$\lim_{x\to 0+}x\psi(x)=-1,\quad\lim_{x\to\infty}\frac{\psi(x)}{\log x}=1.$$
These and many other properties of the $\psi$-function are given, for example, in \cite{1}, \cite{2}.\\

It is our aims in the next theorem to give a refinement of Alzer inequality. We prove,
 \begin{theorem}\label{psi} For all positive real numbers $x\neq 1$,
 $$-\gamma<-\gamma H(x,1/x)<\frac{\gamma^2}{\psi\big(H(x,1/x)\big)}<\psi\Big(1/H(x,1/x)\Big)<H\Big(\psi(x),\psi(1/x)\Big).$$
 Equalities hold if and only if $x=1$
 \end{theorem}
To prove the theorem we need the following lemma.
\begin{lemma} For $x>0$, let
$$\theta(x)=\psi\Big(1/H(x,1/x)\Big).$$
Then $\theta$ admits a unique zero $x_1\in(0,1)$ and $x_1\in(0,1/x_0)$.
\end{lemma}
 {\bf Proof.} For $x>0$,
 $$\theta'(x)=\frac{x^2-1}{2x^2}\psi'\Big(1/H(x,1/x)\Big).$$
 Since $\psi'(x)>0$ for $x>0$, then $\theta$ is strictly decreasing on $(0,1)$ and strictly increasing on $(1,+\infty)$, and $\lim_{x\to 0+}\theta(x)=+\infty$, $\theta(1)=-\gamma$, then $\theta$ admits a unique zero $x_1\in(0,1)$.

 Moreover, for $x\in(0,1)$, $H(x,1/x)\in(0,1)$. Since, $\theta(1)=-\gamma$ and $\theta(2)\simeq -0.227<0$ and $x_0<2$, then $\theta(x_0)<0$. Moreover,
 $\theta(1/x_0)=\theta(x_0)<0$. Then, $\theta(1/x_0)<\theta(x_1)$, using the monotony of $\theta$ on $(0,1)$ we get $x_1<1/x_0$.

{\bf Proof of Theorem \ref{psi}.}
We have, $\psi(1)=-\gamma$, then $-\gamma H(1,1)=H\Big(\psi(1),\psi(1)\Big)=-\gamma$.

Therefore, by the symmetry $x\longleftrightarrow1/x$, it remains to prove the list of the strict inequalities on $(0,1)$.
Let us set $$\sigma(x)=H\Big(\psi(x),\psi(1/x)\Big).$$

 Step 1: $x\in(x_1,1/x_0)$. Since, $\theta(x)<0$ and in Lemma 3 \cite{alz4}, it is proved that $\sigma(x)>0$. This gives the desired result.

 Step 2: $x\in(0,x_1)$. In this case $\sigma(x)>0$ and $\theta(x)>0$. To show that $\theta(x)<\sigma(x)$ is equivalent to prove that
 $$\frac1{\psi(x)}+\frac1{\psi(1/x)}<\frac2{\psi(\frac12(x+1/x)}.$$
 Since, $x_1<1/x_0<x_0$ then, $\psi(x)<0$ and $\psi(\frac12(x+1/x)>0$. therefore,
 $$\frac1{\psi(x)}<\frac1{\psi(\frac12(x+1/x)}.$$
 Moreover, for all $x\in(0,1)$, $1/x>\frac12(x+1/x)$ and the function $\psi$ is strictly increasing on $(0,+\infty)$, then
 $\psi(1/x)>\psi(\frac12(x+1/x),$ or
 $$\frac1{\psi(1/x)}<\frac1{\psi(\frac12(x+1/x)}.$$
 Whence, $$\frac1{\psi(x)}+\frac1{\psi(1/x)}<\frac2{\psi(\frac12(x+1/x)},$$ and
 $\theta(x)<\sigma(x)$.

 Step 3: $x\in(1/x_0,1)$.

 Let $\tau(x)=1/\psi(x)$. Differentiate yields,
 $$\tau''(x)=\frac{2(\psi'(x))^2-\psi(x) \psi''(x)}{(\psi(x))^3}.$$
 Since, for $x>0$, $(\psi'(x))^2+\psi''(x)\geq 0$, then $2(\psi'(x))^2-\psi(x) \psi''(x)\geq -\psi''(x)(2+\psi(x))$. Moreover, for $x>0$, $\psi''(x)<0$, and for $x\geq 1/2$, $\psi(x)>\psi(1/2)>-1.97$ then $2+\psi(x)>0$ and $\psi(x)<0$ for $x\in(1/x_0,x_0)$ therefore, $\tau''(x)<0$ and $\tau$ is strictly concave on $(x_0,1/x_0)$. For all $x\in(1/x_0,x_0)$ we get
 $$\frac12\Big(\frac1{\psi(x)}+\frac1{\psi(1/x)}\Big)<\psi(\frac12(x+1/x))$$

In \cite{alz3}, it is proved that for $y>0$, $y\neq 1$, $\psi(y)\psi(1/y)<\gamma^2$ (Proposition 4). Put $y=H(x,1/x)$, we get
 $$\psi(H(x,1/x))\psi(1/H(x,1/x))<\gamma^2.$$
 Since, for $x\in(0,1)$ we have $H(x,1/x)\in(0,1)$, then
 \begin{equation}\label{40}\psi(H(x,1/x))<0,\end{equation}
and $$\psi(1/H(x,1/x))>\frac{\gamma^2}{\psi(H(x,1/x))}.$$

From equation \eqref{40}, the inequality
$$-\gamma H(x,1/x)<\frac{\gamma^2}{\psi\big(H(x,1/x)\big)}$$
is equivalent to
$$\frac\gamma {H(x,1/x)}+\psi\big(H(x,1/x)\big)<0.$$
For $y>0$ let $u(y)=\gamma/y+\psi(y)$. Then, $u'(y)=\psi'(y)-\gamma/y^2$, since for $y>0$, $\psi'(y+1)=\psi'(y)-1/y^2$, whence
$u'(y)=\psi'(y+1)+(1-\gamma)/y^2>0$. Then, $u$ is strictly increasing and since $u(1)=0$, thus $\gamma/y+\psi(y)<0$ for all $y\in(0,1)$. Moreover, for $x>0$ and $x\neq 1$ we have $H(x,1/x)\in(0,1)$ then $\gamma/H(x,1/x)+\psi(H(x,1/x))<0$.

\section{Zeta function}

For $s>1$, the Riemann zeta function is defined by the infinite series
$$\zeta(s)=\sum_{n=1}^\infty\frac1{n^s}.$$
The function $\zeta$ admits an analytic continuation to the half complex plane $\Re s>0, s\neq 1$. This result is due to Riemann, where he proves that
$$\zeta(s)=\frac{\eta(s)}{1-2^{1-s}},\quad\eta(s)=\sum_{n=1}^\infty\frac{(-1)^{n+1}}{n^s}.$$

On can extends the zeta function to a meromorphic function to all the complex plane in the following way: Since for $\Re z>1$, we have
$$\zeta(z)=\frac1{\Gamma(z)}\int_0^\infty \frac{t^{z-1}}{e^t-1}dt,$$
Splitting the integral
$$\int_0^\infty \frac{t^{z-1}}{e^t-1}dt=\int_0^1 \frac{t^{z-1}}{e^t-1}dt+\int_1^\infty \frac{t^{z-1}}{e^t-1}dt,$$
the second integral converges for all $z\in\Bbb C$ and therefore defines an
entire function which we denote by $F$. On the other hand, the function $z\mapsto 1/(e^z-1)$ has a simple pole at $0$ with residue
$1$. Therefore,
$$\frac1{e^z-1}=\frac1z+G(z),$$
where $G$ is an meromorphic function with first order poles at $2m\pi i$ with integer $m\neq 0$. It follows that for $|z|<2\pi$
$$\frac1{e^z-1}=\frac1z+\sum_{n=0}^\infty c_n z^n.$$
 By Cauchy's estimates, if we fix $0<r<2\pi$ it follows that $|c_n|\leq M/r^n$ for some $M>0$. In particular, there exists $C>0$ such that
 $|c_n|\leq M/2^n$ for $n\in\Bbb N$.

 This implies that, for $\Re z>1$, since the above series converges uniformly on $[0,1]$, we have
 $$\int_0^1 \frac{t^{z-1}}{e^t-1}dt=\frac1{z-1}+\sum_{n=0}^\infty\frac{c_n}{z+n}.$$
 Therefore, we have
 $$\zeta(z) \Gamma(z)=\frac1{z-1}+\sum_{n=0}^\infty\frac{c_n}{z+n}+F(z).$$
 Hence, $z\mapsto \zeta(z)\Gamma(z)$ extends to a meromorphic function in the complex plane with simple pole at 1 and either simple poles or
removable singularities at $z=0,-1,-2,...$ depending if $c_n\neq 0$ or $c_n=0$. Since $\Gamma(z)$ has no zeros, then $z\mapsto1/\Gamma(z)$ is an entire function. This implies that the zeta function $\zeta$ extends to a meromorphic function in the complex plane. Moreover, $\Gamma(1)=1$, then $\zeta$ has a simple pole at $z=1$ with residue $1$. Furthermore, $\Gamma$ has simple pole at $z=0,-1,-2,...$, the function $1/\Gamma$ has simple zero there. It follows that $\zeta$ has removable singularities at $z=0,-1,-2,...$.

  It can therefore be expanded as a Laurent series about $z=1$; the series development is then for $|z-1|<3$
  \begin{equation}\label{1}\zeta(z)=\frac1{z-1}+\sum_{k=0}^\infty\frac{\gamma_n}{n!}(1-z)^n,\end{equation}
  where, for $n\geq 1$ the Stieltjes constants $\gamma_n$ are given by
  $$\gamma_n=\lim_{m\to\infty}\Big(\sum_{k=1}^m\frac{(\log k)^n}{k}-\frac{(\log m)^{n+1}}{n+1}\Big),$$
  and $\gamma_0=\gamma$ is the Euler-Mascheroni constant.

  Many bounds are given for the Stieltjes constants, here is a bound that keeps our interest and will be used later, see for instance \cite{bru}.
  $$|\gamma_{2n}|\leq\frac{4(2n-1)!}{\pi^{2n}},\quad |\gamma_{2n+1}|\leq\frac{2(2n)!}{\pi^{2n+1}},\;\;n=0,1,...$$

The zeta function plays an important role in several branches of mathematics and related areas. In fact, it has interesting applications in the theory of
special functions, in the theory of infinite series, in statistics, in physics and it is
also subject of number theoretic investigation.

It is our main goal to prove a harmonic mean inequality for the zeta function.
 % \begin{theorem} The function $\zeta(x)$ is strictly logarithmic convex on $(0,1)$.  \end{theorem}
\begin{theorem}\label{t1} The function $x\mapsto (\zeta(x)+\zeta(1/x))/(\zeta(x)\zeta(1/x))$ is strictly increasing on $(0,1)$ and strictly decreasing on $(1,+\infty)$ and for all $x>0$
$$H(\zeta(x),\zeta(1/x))<-2.$$
The sign of equality hold if and only if $x = 0$.
\end{theorem}
We use the following result due to Alzer, see Theorem 3.1 and Remark 3.2.
  \begin{proposition}\label{p2} The function $\eta(s)=(1-2^{1-s})\zeta(s)$ is strictly concave on $(0,+\infty)$. Moreover,
  $\eta'(s)>0$ and $1/2<\eta(s)<1$ for all $s>0$.
  \end{proposition}
  As consequence
  \begin{remark}\label{r2} For $s\in(0,1)$, $\zeta(s)<0$ and $\zeta'(s)<0$
  \end{remark}
  Among others results states in this work we prove.
\begin{proposition}\label{pr6} \
\begin{enumerate}
\item [$(1)$] The function $\zeta(x)$ is strictly completely monotonic on $(1,\infty)$ and the function $-\zeta(x)$ is absolutely monotonic on $(0,1)$.
\item [$(2)$] The function $g(x)=x\zeta'(x)-(1/x)\zeta'(1/x)$ is continuous at $x=1$ with $g(1)=0$, $g(x)<0$ for $x\in(0,1)$ and $g(x)>0$ on $(1,+\infty)$.
\end{enumerate}
\end{proposition}
As consequence. The function $\zeta$ is strictly concave on $(0,1)$ and strictly convex on $(1,+\infty)$.

To prove these results we need the following lemmas
\begin{lemma}\label{lem01}
For all $x\geq 2$ and all $a\geq 2$,
$$4^{1-x} ((x-1) \log (4)+1)\leq \frac x2,$$ and
$$\frac{\log a}{a^x}\leq \frac1{x^2}.$$
\end{lemma}
{\bf Proof.} Differentiate the function $T(x)=4^{1-x} ((x-1) \log (4)+1)-x/2$ yields $T'(x)=-2^{-2 x-1}(4^x+8\log ^2(4)(x-1))<0$. Then $T$ is strictly decreasing and $T(2)=1/4(\log 4+1)-1<0.$

Differentiate the function $S(x)=x^2\log a/a^x-1$ yields $S'(x)=xa^{-x} \log (a) (2-x \log (a))$.
Therefore, $S(x)\leq S(2/\log a)=(4/\log a)e^{-2}-1\leq0$ if and only if $a\geq e^{4/e^2}\simeq 1.71$. This completes the proof of the lemma.
\begin{lemma}\label{lem0} \
\begin{enumerate}
\item [$(1)$] For $x\in(0,1)$,
$$\zeta'(x)< -\frac1{(1-x)^2}-\gamma_1-\gamma_2(1-x).$$
\item [$(2)$] For all $x\in(1,2]$, $$\zeta'(x)>-\frac{1}{(x-1)^2}+\gamma _2 (x-1)-\gamma _1-\frac{1}{2} \gamma _3 (x-1)^2.$$
\item [$(3)$] For $x\geq 1$,
$$-\zeta'(x)\leq \frac{4^{1-x} ((x-1) \log (4)+1)}{(x-1)^2}+\frac{\log 2}{2^x}+\frac{\log 3}{3^x}+\frac{\log 4}{4^x}.$$
%3) For every $x\in(1,4)$,
%$$\zeta'(x)>-\frac{1}{(x-1)^2}+\gamma _2 (x-1)-\gamma _1-\frac{1}{2} \gamma _3 (x-1)^2.$$
\end{enumerate}
\end{lemma}
{\bf Proof}
1) Let $$\theta(x)=\zeta'(x)+\frac1{(1-x)^2}+\gamma_1+\gamma_2(1-x).$$ By successive differentiation we get
$$\theta'(x)=-\gamma _2+\zeta ''(x)+\frac{2}{(1-x)^3},$$
$$\theta''(x)=\zeta '''(x)+\frac{6}{(1-x)^4},$$
therefore,
$$\theta''(x)=-\sum_{k=3}^\infty\frac{k(k-1)(k-2)\gamma_k}{k!}(1-x)^{k-3},$$
then,
$$\theta''(x)=-\sum_{k=3}^{2n}\frac{k(k-1)(k-2)\gamma_k}{k!}(1-x)^{k-3}+R_{2n}(x),$$
where, $$R_{2n}(x)=-\sum_{k=2n+1}^{\infty}\frac{k(k-1)(k-2)\gamma_k}{k!}(1-x)^{k-3}.$$
Since, $\gamma_{2k}\leq 4(2k-1)!/\pi^{2k}$ and $\gamma_{2k+1}\leq 4(2k)!/\pi^{2k+1}$. Then, for $x\in(0,1)$
$$\begin{aligned}R_{2n}(x)&=-\sum_{k=2n+1}^{\infty}\frac{(k-1)(k-2)\gamma_k}{(k-1)!}(1-x)^{k-3}\\&\leq 4\sum_{k=n+1}^{\infty}\frac{(2k-1)(2k-2)}{\pi^{2k}}(1-x)^{2k-3}+
2\sum_{k=n}^{\infty}\frac{(2k)(2k-1)}{\pi^{2k+1}}(1-x)^{2k-2}.\end{aligned}$$
Hence,
$$\begin{aligned}R_{2n}(x)&\leq 4 \pi ^{-2 n-1} \frac{(a_n x^4+b_n x^3+c_n x^2+d_nx+e_n) (1-x)^{2 n-2}}{\left(\pi ^2-(x-1)^2\right)^3}\\
&+8 \pi ^{-2 n-2}\frac{(a'_n x^6+b'_n x^5+c'_n x^4+d'_nx^3+e'_nx^2+f'_nx+g'_n) (1-x)^{2 n-3}}{\left(\pi ^2-(x-1)^2\right)^3}\end{aligned}$$
$$a_n=2 \pi ^2 n^2-5 \pi ^2 n+3 \pi ^2,b_n=-8 \pi ^2 n^2+20 \pi ^2 n-12 \pi ^2,$$
$$c_n=-4 \pi ^4 n^2+12 \pi ^2 n^2+6 \pi ^4 n-30 \pi ^2 n+\pi ^4+18 \pi ^2,d_n=8 \pi ^4 n^2-8 \pi ^2 n^2-12 \pi ^4 n+20 \pi ^2 n-2 \pi ^4-12 \pi ^2 $$
$$e_n=2 \pi ^6 n^2-4 \pi ^4 n^2+2 \pi ^2 n^2-\pi ^6 n+6 \pi ^4 n-5 \pi ^2 n+\pi ^4+3 \pi ^2$$

$$a'_n=2 \pi ^2 n^2-3 \pi ^2 n+\pi ^2,b'_n=-12 \pi ^2 n^2+18 \pi ^2 n-6 \pi ^2$$
$$c'_n=-4 \pi ^4 n^2+30 \pi ^2 n^2+2 \pi ^4 n-45 \pi ^2 n+3 \pi ^4+15 \pi ^2 ,$$
$$d'_n=16 \pi ^4 n^2-40 \pi ^2 n^2-8 \pi ^4 n+60 \pi ^2 n-12 \pi ^4-20 \pi ^2 $$
$$e'_n=2 \pi ^6 n^2-24 \pi ^4 n^2+30 \pi ^2 n^2+\pi ^6 n+12 \pi ^4 n-45 \pi ^2 n+18 \pi ^4+15 \pi ^2$$
$$f'_n=-4 \pi ^6 n^2+16 \pi ^4 n^2-12 \pi ^2 n^2-2 \pi ^6 n-8 \pi ^4 n+18 \pi ^2 n-12 \pi ^4-6 \pi ^2$$
$$g'_n=2 \pi ^6 n^2-4 \pi ^4 n^2+2 \pi ^2 n^2+\pi ^6 n+2 \pi ^4 n-3 \pi ^2 n+3 \pi ^4$$
For $n=5$, we get
$\pi^9(\pi ^2-(x-1)^2)^3 \theta''(x)\leq P(x),$
where,
$$\begin{aligned}P(x)&=(1-x)^{13} \left(\frac{\gamma _{10}}{5040}+288\right)+(1-x)^{12} \left(\frac{\gamma _9}{720}+112 \pi \right)+(1-x)^{11} \left(\frac{\gamma _8}{120}-\frac{\pi ^2 \gamma _{10}}{1680}-696 \pi ^2\right)\\&+(1-x)^{10} \left(\frac{\gamma _7}{24}-\frac{\pi ^2 \gamma _9}{240}-276 \pi ^3\right)+(1-x)^9 \left(\frac{\gamma _6}{6}-\frac{\pi ^2 \gamma _8}{40}+\frac{\pi ^4 \gamma _{10}}{1680}+440 \pi ^4\right)\\&+(1-x)^8 \left(\frac{\gamma _5}{2}-\frac{\pi ^2 \gamma _7}{8}+\frac{\pi ^4 \gamma _9}{240}+180 \pi ^5\right)+(1-x)^7 \left(\gamma _4-\frac{\pi ^2 \gamma _6}{2}+\frac{\pi ^4 \gamma _8}{40}-\frac{\pi ^6 \gamma _{10}}{5040}\right)\\&+(1-x)^6 \left(\gamma _3-\frac{3 \pi ^2 \gamma _5}{2}+\frac{\pi ^4 \gamma _7}{8}-\frac{\pi ^6 \gamma _9}{720}\right)+(1-x)^5 \left(-3 \pi ^2 \gamma _4+\frac{\pi ^4 \gamma _6}{2}-\frac{\pi ^6 \gamma _8}{120}\right)\\&+(1-x)^4 \left(-3 \pi ^2 \gamma _3+\frac{3 \pi ^4 \gamma _5}{2}-\frac{\pi ^6 \gamma _7}{24}\right)+(1-x)^3 \left(3 \pi ^4 \gamma _4-\frac{\pi ^6 \gamma _6}{6}\right)\\
&+(1-x)^2 \left(3 \pi ^4 \gamma _3-\frac{\pi ^6 \gamma _5}{2}\right)-\pi ^6 (1-x) \gamma _4-\pi ^6 \gamma _3.\end{aligned}$$

An application of Sturm's theorem reveals that the polynomial $P(x)$ has no zero on $(0,1)$ and since $P(1)=-\pi^6\gamma_3<0$, therefore $\theta''(x)<0$ for all $x\in(0,1)$. So, $\theta'$ is strictly decreasing and $\theta'(1)=0$ then $\theta$ is strictly increasing on $(0,1)$ and $\theta(1)=0$. Therefore,
$\theta(x)<0$ and the result follows.

2) Let $$\tau(x)=\zeta'(x)+\frac{1}{(x-1)^2}-\gamma _2 (x-1)+\gamma _1+\frac{1}{2} \gamma _3 (x-1)^2.$$
Then,
$$\tau''(x)=\zeta'''(x)+\frac{6}{(x-1)^4}+\gamma_3=\theta''(x)+\gamma_3.$$
Then,
$$\tau''(x)=-\sum_{k=4}^{2n}\frac{k(k-1)(k-2)\gamma_k}{k!}(1-x)^{k-3}+R_{2n}(x),$$
where, $$R_{2n}(x)=-\sum_{k=2n+1}^{\infty}\frac{k(k-1)(k-2)\gamma_k}{k!}(1-x)^{k-3}.$$
Since, $\gamma_{2k}\geq -4(2k-1)!/\pi^{2k}$ and $\gamma_{2k+1}\leq -4(2k)!/\pi^{2k+1}$. Then, for $x\geq 1$
$$\begin{aligned}R_{2n}(x)&=-\sum_{k=2n+1}^{\infty}\frac{(k-1)(k-2)\gamma_k}{(k-1)!}(1-x)^{k-3}\\&\geq -4\sum_{k=n+1}^{\infty}\frac{(2k-1)(2k-2)}{\pi^{2k}}(x-1)^{2k-3}-
2\sum_{k=n}^{\infty}\frac{(2k)(2k-1)}{\pi^{2k+1}}(x-1)^{2k-2}.\end{aligned}$$
For $n=5$, we get
$\pi^{10}(\pi ^2-(x-1)^2)^3 \tau''(x)\geq (x-1)Q(x),$
where,
$$\begin{aligned}Q(x)&=\left(-\frac{\pi ^{10} \gamma _{10}}{5040}-288\right) (x-1)^{13}+\left(\frac{\pi ^{10} \gamma _9}{720}-112 \pi \right) (x-1)^{11}\\
&+\left(-\frac{\pi ^{10} \gamma _8}{120}+\frac{\pi ^{12} \gamma _{10}}{1680}+696 \pi ^2\right) (x-1)^{10}+\left(\frac{\pi ^{10} \gamma _7}{24}-\frac{\pi ^{12} \gamma _9}{240}+276 \pi ^3\right) (x-1)^{9}\\
&+\left(-\frac{\pi ^{10} \gamma _6}{6}+\frac{\pi ^{12} \gamma _8}{40}-\frac{\pi ^{14} \gamma _{10}}{1680}-440 \pi ^4\right) (x-1)^8
+\left(\frac{\pi ^{10} \gamma _5}{2}-\frac{\pi ^{12} \gamma _7}{8}+\frac{\pi ^{14} \gamma _9}{240}-180 \pi ^5\right) (x-1)^7\\
&+\left(-\pi ^{10} \gamma _4+\frac{\pi ^{12} \gamma _6}{2}-\frac{\pi ^{14} \gamma _8}{40}+\frac{\pi ^{16} \gamma _{10}}{5040}\right) (x-1)^6+\left(\frac{1}{2} (-3) \pi ^{12} \gamma _5+\frac{\pi ^{14} \gamma _7}{8}-\frac{\pi ^{16} \gamma _9}{720}\right) (x-1)^5\\&
+\left(3 \pi ^{12} \gamma _4-\frac{\pi ^{14} \gamma _6}{2}+\frac{\pi ^{16} \gamma _8}{120}\right) (x-1)^4+\left(\frac{3 \pi ^{14} \gamma _5}{2}-\frac{\pi ^{16} \gamma _7}{24}\right) (x-1)^3\\
&+\left(\frac{\pi ^{16} \gamma _6}{6}-3 \pi ^{14} \gamma _4\right) (x-1)^2-\frac{1}{2} \pi ^{16} \gamma _5 (x-1)+\pi ^{16} \gamma _4 .\end{aligned}$$
An application of Sturm's theorem reveals that the polynomial $Q(x)$ has no zero on $(1,2)$ and since $Q(1)=\pi^{16}\gamma_4>0$, therefore $\tau''(x)>0$ for all $x\in(1,2)$. So, $\tau'$ is strictly increasing and $\tau'(1)=0$ then $\tau$ is strictly increasing on $(1,2)$ and $\tau(1)=0$. Therefore,
$\tau(x)>0$ and the result follows.

3) For every $x\geq 1$ the function $t\mapsto \log t/t^x$ is strictly decreasing on $[e,+\infty)$. Therefore, for all $k\geq 3$
$$\frac{\log (k+1)}{(k+1)^x}\leq \int_k^{k+1}\frac{\log t}{t^x} dt\leq \frac{\log k}{k^x},$$
Then $$-\zeta'(x)\leq \int_3^\infty \frac{\log t}{t^x}  dt+\frac{\log 2}{2^x}+\frac{\log 3}{3^x}=\frac{3^{1-x} ((x-1) \log (3)+1)}{(x-1)^2}+\frac{\log 2}{2^x}+\frac{\log 3}{3^x}.$$
This completes the proof of the lemma.

As consequence, we have the following inequalities
\begin{corollary}\label{cor1} For all $x\in(0,1)$
$$\zeta''(x)>\gamma_2-\frac2{(1-x)^3}.$$
$$\zeta(x)>\frac1{x-1}+\gamma+\gamma_1 (1-x)+\frac{\gamma_2}2(1-x)^2,$$
and
$$\zeta(x)<\frac1{x-1}-\frac{2\gamma_1+\gamma_2-1}2+\gamma_1 (1-x)+\frac{\gamma_2}2(1-x)^2,$$
\end{corollary}
{\bf Proof of the proposition \ref{pr6}.}

1)  Recall that for
$x>1$, $$\zeta(x)=\sum_{k=1}^\infty\frac1{k^x}.$$
By termwise differentiation of the series we get $$(-1)^n\zeta^{(n)}(x)=\sum_{k=1}^\infty\frac{(\log k)^n}{k^x}>0.$$

 By termwise differentiation of the Laurent series \eqref{1}
 $$(-1)^n\zeta^{(n)}(x)=\frac{n!}{(x-1)^{n+1}}+\sum_{k=n}^\infty\frac{k(k-1)...(k-n+1)\gamma_k}{k!}(1-x)^{k-n}.$$
 Using the inequality of Lavrik  $|\gamma_k|\leq k!/2^{k+1}$ for $k\geq 1$, we get for $x\in(0,1)$
 $$\big|\sum_{k=n}^\infty\frac{k(k-1)...(k-n+1)\gamma_k}{k!}(1-x)^{k-n}\big|\leq \frac1{2^{n+1}}\sum_{k=n}^\infty k(k-1)...(k-n+1)(\frac{1-x}2)^{k-n}.$$
 The sum in the right hand side is nothing rather than the $n$-derivative of the function
 $((-1)^n/2)\sum_{k=0}^\infty(\frac{1-x}2)^{k}$. Therefore, for $x\in(0,1)$
 $$\big|\sum_{k=n}^\infty\frac{k(k-1)...(k-n+1)\gamma_k}{k!}(1-x)^{k-n}\big|\leq \frac{n!}{(1+x)^{n+1}}.$$
 Then,
 \begin{equation}\label{e3} -\frac{(-1)^nn!}{(1-x)^{n+1}}-\frac{n!}{(1+x)^{n+1}}\leq (-1)^n\zeta^{(n)}(x)\leq \frac{n!}{(1+x)^{n+1}}-\frac{n!(-1)^n}{(1-x)^{n+1}}.\end{equation}
 Whence,  \begin{equation}\label{e2} -\frac{n!}{(1-x)^{n+1}}(1+(-1)^{n})\leq (-1)^n\zeta^{(n)}(x)\leq \frac{n!}{(1-x)^{n+1}}(1+(-1)^{n+1}).\end{equation}
 Then $\zeta^{(n)}(x)\leq 0$ for all $n$ and all $x\in(0,1)$. Which completes the proof of the lemma.

2) Recall that $g(x)=x\zeta'(x)-\frac1x\zeta'(\frac1x),$
then by Lemma \ref{lem0}, we get for $x\in(0,1)$,
$$g(x)\leq-\frac1x\zeta'(\frac1x) -\frac x{(1-x)^2}-\gamma_1 x-\gamma_2x(1-x):=f(\frac1x).$$

For $x>1$
$$f(x)=-x\zeta'(x) -\frac{x}{(1-x)^2}-\frac{\gamma_1} x-\gamma_2\frac{x-1}{x^2}.$$
Assume $x\in(2,\infty)$, then by Lemma \ref{lem0}, we get,
$$f(x)\leq \frac{4^{1-x} ((x-1) \log (4)+1)}{(x-1)^2}+\frac{\log 2}{2^x}+\frac{\log 3}{3^x}+\frac{\log 4}{4^x}-\frac{x}{(1-x)^2}-\frac{\gamma_1}x-\gamma_2\frac{x-1}{x^2}.$$
Moreover, From Lemma \ref{lem01}, $4^{1-x} ((x-1) \log (4)+1)\leq x/2$ and $2^{-x} \log 2 +3^{-x} \log 3 +4^{-x}\log 4 \leq 3/x^2$. Therefore, for $x\geq 2$
$$f(x)\leq\frac3{x^2}-\frac{x}{2(1-x)^2}-\frac{\gamma_1}x-\gamma_2\frac{x-1}{x^2}.$$
%$$x^2(x-1)^2f(\frac1x)\leq3(x-1)^2-\frac{x^3}{2}-\gamma_1x(x-1)^2-\gamma_2(x-1)^3.$$
Then, $$x^2(x-1)^2f(x)\leq\left(-\gamma _1-\gamma _2-\frac{1}{2}\right) x^3+\left(2 \gamma _1+3 \gamma _2+3\right) x^2+\left(-\gamma _1-3 \gamma _2-6\right) x+\gamma _2+3:=P_1(x).$$
Differentiate yields
$$P'_1(x)=-\frac{1}{2} 3 \left(2 \gamma _1+2 \gamma _2+1\right) x^2+\left(4 \gamma _1+6 \gamma _2+6\right) x-\gamma _1-3 \left(\gamma _2+2\right).$$
and $P'_1$ admits a unique root $x_1$ in $[2,+\infty)$ and $P'_1(2)\simeq0.393>0$. Therefore, $P_1$ is strictly increasing on $(2,x_1)$ and strictly decreasing on $(x_1,+\infty)$. Moreover, by an easy mathematica software computation we have $P_1(x_1)\simeq-0.535<0$. Thus, for all $x\geq 2$, $fx)<0$ and $g(x)<0$ for all $x\in(0,1/2]$.

Assume $x\in(1,2]$.  By using Lemma \ref{lem0},
$$f(x)\leq-\gamma _2 x(x-1)+\gamma _1x+\frac{1}{2} \gamma _3 x(x-1)^2-\frac{\gamma_1} x-\gamma_2\frac{x-1}{x^2}.$$
Then $$\frac{2x^2}{x-1}f(x)\leq \left(2 \gamma _2+\gamma _3\right) x^3+\left(2 \gamma _1-\gamma _3\right) x^2+2 \gamma _1 x-2 \gamma _2:=v(x).$$

Differentiate yields
$$v'(x)=3 \left(2 \gamma _2+\gamma _3\right) x^2+\left(4 \gamma _1-2 \gamma _3\right) x+2 \gamma _1.$$
Numerical computation show that $v'(x)$ admits two negative roots and since, $v'(1)\simeq-0.492$, then $v'(x)<0$ on $(1,2]$ and $v(x)$ is strictly decreasing on $(1,2)$ and $v(1)\simeq -0.291$, thus $v(x)<0$ on $(1,2]$ or equivalently $f(x)<0$ on $(1,2]$ which implies that $g(x)<0$ on $[1/2,1)$.
Moreover, $g(1/x)=-g(x)$ then $g(x)>0$ for all $x>1$. The continuity of $g$ at $x=1$ follows from equation \eqref{1}, indeed for $|x-1|<3$
$$g(x)=\sum_{k=1}^\infty\frac{\gamma_k}{(k-1)!}(1-x)^{k-1}-\sum_{k=1}^\infty\frac{\gamma_k}{(k-1)! x^k}(x-1)^{k-1}=\gamma_1-\frac{\gamma_1}x+o(x-1).$$
This completes the proof of the proposition.

\begin{lemma}\label{lem1} The functions $|\zeta(x)|$ is strictly logarithmic convex on $(0,1)$.% $u(x)=x\zeta'(x)$ and $v(x)=x\zeta'(x)\zeta^2(1/x)$ are strictly increasing on $(1,+\infty)$.

For all $x\in(0,1)$
$$\zeta(x)\zeta(1-x)>(\zeta(1/2))^2.$$
Equality holds if and only if $x=1/2$.
\end{lemma}
{\bf Proof.}
For every $x\in(0,1)$,
$$(\log|\zeta(x)|)''=\frac{\zeta''(x)\zeta(x)-(\zeta'(x))^2}{(\zeta(x))^2}.$$
Form equation \eqref{e2}, we have
$\zeta''(x)<2/(1+x)^3-2/(1-x)^3$, $\zeta(x)<1/(1+x)-1/(1-x)$ and $0>\zeta'(x)>-1/(1+x)^2-1/(1-x)^2$. Therefore,
$$\zeta''(x)\zeta(x)-(\zeta'(x))^2\geq \frac{4(x^4+4 x^2-1)}{(x^2-1)^4}.$$
The unique root of the polynomial $x^4+4 x^2-1$ in the interval $(0,1)$ is $\sqrt5-2$ and $x^4+4 x^2-1>0$ on $(\sqrt5-2,1)$.

Furthermore, From Corollary \ref{cor1} and by integration we get for $x\in(0,1)$
$$\zeta'(x)>-\frac1{(1-x)^2}+\gamma_2 x+1-\frac{\log(2\pi)}2.$$
Since $\gamma_2 x+1-\frac{\log(2\pi)}2>\gamma_2 +1-\frac{\log(2\pi)}2\simeq 0.071$ then
$$\zeta'(x)>-\frac1{(1-x)^2}.$$
So, from corollary \ref{cor1}, we get
$$\zeta''(x)\zeta(x)-(\zeta'(x))^2>(\frac1{(1-x)^2}-\gamma_2)(\frac1{1-x}+\frac1{1+x})-\frac1{(1-x)^4}.$$
Since, $\gamma_2<0$, then

$$\zeta''(x)\zeta(x)-(\zeta'(x))^2>\frac{1-3x}{(1-x)^4(1+x)}.$$
Then, $\zeta''(x)\zeta(x)-(\zeta'(x))^2>0$ for $x\in(0,1/3)$. This completes the proof of the lemma.

Let $\varphi(x)=\log|\zeta(x)|+\log|\zeta(1-x)|$. The function $\varphi$ is convex, then the function $\varphi'(x)=\frac{\zeta'(x)}{\zeta(x)}-\frac{\zeta'(1-x)}{\zeta(1-x)}$ is strictly increasing, moreover, $\varphi'(1/2)=0$, then $\varphi$ is strictly decreasing on $(0,1/2)$ and strictly increasing on $(1/2,1)$ and
$\varphi(1/2)=2\log|\zeta(1/2)|$. So for $x\in(0,1)$, $x\neq 1/2$, $$\zeta(x)\zeta(1-x)>(\zeta(1/2))^2.$$
\begin{proposition} \
\begin{enumerate}
\item The function $1/\zeta$ is strictly concave on $(0,1)$.

\item For all $x\in(0,1)$
$$\zeta(1/2)<H(\zeta(x),\zeta(1-x))<-1.$$
Equality holds if and only if $x=0,1$ and $x=1/2$.
\item For all $x\in(0,1)$
$$\frac{\log 4}{1+\log 4}<H(\eta(x),\eta(1-x))<(1-\sqrt 2)\zeta(1/2).$$
Equality holds if and only if $x=0,1$ and $x=1/2$.
\end{enumerate}
\end{proposition}
{\bf Proof.}  a) For $x\in(0,1)$ we have
\begin{equation}\label{co}\Big(\frac1{\zeta(x)}\Big)''=-\frac{\zeta''(x)\zeta(x)-2\zeta''(x)}{\zeta(x)^3}.\end{equation}
From Corollary \ref{cor1} and Lemma \ref{lem0}, we get
$$\zeta''(x)\zeta(x)-2(\zeta'(x))^2<-\frac{2 \gamma}{(1-x)^3}-\frac{6 \gamma _1}{(1-x)^2}-\frac{3\gamma^2 _2}{2} (1-x)^2 -3\gamma _2 \gamma _1(1-x)-\frac{6\gamma _2}{1-x}-2\gamma_1^2+\gamma  \gamma _2:=g(x),$$
say. Differentiate twice, we get
$$-\frac13(1-x)^5g''(x)=(1-x)(\gamma^2_2(1-x)^4+4  \gamma _2(1-x)+12 \gamma _1)+8 \gamma:=(1-x)h(x)+8 \gamma,$$
say. Differentiate yields,
$h'(x)=-4(\gamma^2_2(1-x)^3+\gamma _2)$. It easy to see that $h'$ is strictly increasing and $h'(0)=-4\gamma_2(\gamma_2+1)>0$. Therefore, $h$ is strictly increasing and $h(1)=12\gamma_1<0$. Therefore, the function $-\frac13(1-x)^5g''(x)$ is strictly increasing and $-\frac13g''(0)=\gamma^2_2+4\gamma_2+12\gamma_1+8\gamma\simeq 3.70527$. Then, $g''(x)<0$ and $g'(x)$ is strictly decreasing. Since,
$$g'(x)=3 (\gamma^2_2 (1-x)-\frac{2 \gamma _2}{(1-x)^2}-\frac{4 \gamma _1}{(1-x)^3}+\gamma _1 \gamma _2-\frac{2 \gamma }{(1-x)^4}),$$
and $g'(0)\simeq -2.52896$. Then, $g$ is strictly decreasing on $(0,1)$ and $g(0)\simeq -0.677849$. Whence, $\zeta''(x)\zeta(x)-2(\zeta'(x))^2<0$ on $(0,1)$, moreover, $\zeta(x)<0$ on $(0,1)$, by equation \eqref{co}, we deduce that the function $1/\zeta$ is strictly concave.

b) For $x\in(0,1)$, let $\varphi(x)=1/\zeta(x)+1/\zeta(1-x)$. Differentiate yields,
$$\varphi'(x)=-\frac{\zeta'(x)}{\zeta(x)^2}+\frac{\zeta'(1-x)}{\zeta(1-x)^2}.$$
By the first item we saw that the function $-\zeta'(x)/(\zeta(x)^2)$ is strictly decreasing on $(0,1)$. So, for $x\in(0,1/2)$, $x<1-x$ we get
$\varphi'(x)<0$ then $\varphi$ is strictly decreasing on $(0,1/2)$ and by the symmetry $\varphi(1-x)=\varphi(x)$ the function $\varphi$ is strictly increasing on $(1/2,1)$. Whence, $\varphi(1/2)=2/\zeta(1/2)<\varphi(x)<\varphi(0)=-2$ with equality if and only if $x=0,1$ and $x=1/2$.

b) From Proposition \ref{p2}, the function $\eta$ is strictly concave hence $-\log(\eta)$ is strictly convex. Since, the exponential function is increasing and convex, then $1/\eta$ is strictly convex too. The rest of the proof follows the same line as b).
%Differentiate yields
%$$u'(x)=\zeta'(x)+x\zeta''(x)=\sum_{n=2}^\infty\frac{(x\log n-1)\log n}{n^x}. $$
%Firstly, we have $x\log n-1>0$ for all $x\geq 1$ and $n\geq 3$.
%Then, $$u'(x)\geq\frac{\log 2}{2^x}(x\log 2-1+\frac{x\log 4-1}{2^{x-1}}).$$
%Let $f(x)=x\log 2-1+(x\log 4-1)/2^{x-1}$ then $f'(x)=2^{-x} \log (2)(2^x-4 x \log (2)+6)=2^{-x} \log (2) g(x),$ and $g'(x)=(2^{x-2}-1)\log 4$. Then, $g$ is strictly decreasing on $(1,2)$ and strictly increasing on $(2,+\infty)$ and $g(2)=10-8\log (2)\simeq4.45 >0$. Therefore, $f$ is strictly increasing on $(1,+\infty)$ and $f(1)=3\log2-2\simeq0.079>0$. It follows that $f(x)>0$ and $u'(x)>0$. This completes the proof of the lemma.
%We saw that $u(x)$ is strictly increasing and negative on $(1,+\infty)$, then $-u(x)=-x\zeta'(x)$ is strictly decreasing and positive. Moreover, the function $\zeta(x)$ is strictly decreasing and positive on $(1,+\infty)$, then the function $x\mapsto 1/\zeta^2(1/x)$ is strictly decreasing and positive on $(1,+\infty)$. Therefore,
%$-v(x)=-x\zeta'(x)/\zeta^2(1/x),$ is strictly decreasing and positive on $(1,+\infty)$.  This completes the proof of the lemma.

\begin{proposition}\label{pr4} For $s>0$, $s\neq 1$ the function $$\varphi(s)=\zeta(s)+\zeta(1/s).$$
is strictly decreasing on $(0,1)$ and strictly increasing on $(1,+\infty)$
and  $$2\gamma-1<\zeta(s)+\zeta(1/s)<1/2.$$
The upper and lower bounds are sharp.
\end{proposition}
 {\bf Proof.} Differentiate yields

 $$\varphi'(s)=\frac1{s}(s\zeta'(s)-\frac1s\zeta'(1/s))=\frac{g(s)}{s},$$
 and the result follows from Lemma \ref{lem0}. Moreover, for $s$ close to $1$
  $$\varphi(s)=2\gamma-1+\sum_{n=1}^\infty \frac{\gamma_n}{n!}(1-\frac{(-1)^n}{s^n})(1-s)^n.$$
  Then, $\lim_{s\to 1}\varphi(s)=2\gamma-1$. Moreover, it is well known that $\zeta(0)=-1/2$.
  Putting all this together, we get $2\gamma-1<\varphi(s)<1/2$ for all $s\in(0,+\infty), s\neq 1$.

  \begin{proposition}\label{pr4} \

  \begin{enumerate}
\item  The function $G_1(x)=(x-1)\zeta(x)$ is strictly log-concave on $(0,+\infty)$ and the function $H_{a,b}(x)=x^a(x-1)^b\zeta(x)$ strictly log-concave on $(1,+\infty)$ if and only if $a\geq 0$ and $b\geq 1$.

 \item The function $G_2(x)=x(1-x)\zeta(x)\zeta(1-x)$ is strictly increasing on $(0,1/2)$ and strictly decreasing on $(1/2,1)$.

  For $x\in(0,1)$,  $$\frac{1}{2x(1-x)}<\zeta(x)\zeta(1-x)<\frac{\zeta(1/2)^2}{4x(1-x)}.$$
  The bounds are sharp.

% The function $\frac{(x-1)\zeta(x)}{x^\alpha}$ is strictly decreasing on $(0,1)$ if and only if $\alpha\geq \gamma$.
 \end{enumerate}
%  $G_1(x)=(1-2^{1-x})(1-2^{x})\zeta(x)\zeta(1-x)=\eta(x)\eta(1-x)$ is strictly log-concave on $(0,1)$. Moreover,
%$G_1$ is strictly increasing on $(0,1/2)$ and strictly decreasing on $(1/2,1)$ and
%$$\frac{\log 2}{2(1-2^{1-x})(1-2^{x})}<\zeta(x)\zeta(1-x)<\frac{(3 - 2 \sqrt 2) \zeta(1/2)^2}{(1-2^{1-x})(1-2^{x})}.$$

%The function $G_2(x)=\frac{x(1-x)}{(1-2^{1-x})(1-2^{x})}$ is strictly increasing on $(0,1/2)$ and strictly decreasing on $(1/2,1)$.
\end{proposition}
{\bf Proof}
 1) Let $\varphi(x)=\log((x-1)/(1-2^{1-x}))$. It is not difficult to show that $\varphi$ of class at least two on $(0,+\infty)$ and $\varphi(1)=\log\log2$, $\varphi'(1)=\log2/2$, $\varphi''(1)=(\log 2)^2/12$. Differentiate yields
 $$\varphi'(x)=\frac1{x-1}-\frac{2^{1-x}\log 2}{1-2^{1-x}}=\frac{1-2^{1-x}-(x-1)2^{1-x}\log 2}{(x-1)(1-2^{1-x})}=\frac{u(x)}{v(x)}.$$
 Moreover,
 $$\frac{u'(x)}{v'(x)}=\frac{(x-1)2^{1-x}}{1-2^{1-x}+(x-1)2^{1-x}\log 2}(\log2)^2,$$
and,
  $$\Big(\frac{u'(x)}{v'(x)}\Big)'=2(\log 2)^2\frac{(2^x+2^x \log (2)-2^x x \log (2)-2)}{(2^x+x \log (4)-2-\log (4))^2}.$$
  By easy computation, one see that $x\mapsto 2^x(1+ \log (2)(1-x))-2$ is strictly increasing on $(0,1)$ and strictly decreasing on $(1,+\infty)$ and negative. Therefore, $\frac{u'(x)}{v'(x)}$ is strictly decreasing and since, $u(1)=v(1)=0$ then, the function $u(x)/v(x)$ is strictly decreasing $(0,1)$ and on $(1,+\infty)$ and $\varphi$ is strictly concave on $(0,+\infty)$.
  It is proved by Wang that $\eta$ is logarithmic concave on $(0,\infty)$. Since, the sum of log-concave functions is log-concave. Therefore, the function $(x-1)\eta(x)/(1-2^{1-x})=(x-1)\zeta(x)$ is strictly log-concave on $(0,+\infty)$.

  If $H_{a,b}$ is strictly log-concave, then $(\log H_{a,b})'=\frac ax+\frac b{x-1}+\frac{\zeta'(x)}{\zeta(x)}$ is strictly decreasing. Since
  $\lim_{x\to 1}(1/(x-1)+\zeta'(x)/\zeta(x))=\gamma$ then $b-1\geq 0$. Moreover,  $\lim_{x\to 0}(b/(x-1)+\zeta'(x)/\zeta(x)=-b+\log(2\pi)$, then $a\geq 0$.

  The converse, Assume $a\geq 0$, $b\geq 1$, we have
  $$(\log H_{a,b}(x))''=-\frac a{x^2}-\frac{(b-1)}{(x-1)^2}+(\log(G_1(x))''<0,$$
and $H_{a,b}$ is strictly log-concave.

2)  Let $\psi(x)=\log((x-1)\zeta(x))$, then $\log G_2(x)=\log\psi(x)+\log\psi(1-x)$. Since, $\psi$ is strictly concave on $(0,1)$, then $(\log G_2)'(x)$ is strictly decreasing on $(0,1)$. Moreover, $(\log G_2)'(1/2)=0$. Therefore, $G_2$ is strictly increasing on $(0,1/2)$ and strictly decreasing on $(1/2,1)$. Furthermore, $$G_2(0)=\frac12<G_2(x)<G_2(1/2)=(\frac12\zeta(\frac12))^2.$$
%  , and by Alzer, a stronger result was proved, which states that $\eta$ is strictly concave.
%Let $t\in(0,1)$ and $x\in(0,1)$, then
%$\eta(1-(tx+(1-t)y)=\eta(t(1-x)+(1-t)(1-y))>t\eta(1-x)+(1-t)\eta(1-y)$ therefore, $\eta(1-x)$ is strictly concave and $\log\eta(1-x)$ is strictly concave. Since the sum of strictly concave function is strictly concave. Therefore, $\log\eta(x)+\log\eta(1-x)$ is strictly concave.

%From the log-concavity, we deduce that $(\log G_1)'(x)$ is strictly decreasing and since, $(\log G_1)'(1/2)=0$. Then, $\log G_1$ is strictly increasing on $(0,1/2)$ and strictly decreasing on $(1/2,1)$ and similarly for $G_1$.

%2) Firstly remark that $G_2(1-x)=G_2(x)$. So, we check the monotony on $(0,1/2)$. Let
%$\varphi(x)=\log x-\log(2^x-1)$, then $\varphi'(x)=(2^x-1-2^xx\log 2)/x(2^x-1)=u(x)/v(x)$ and
%$u'(x)/v'(x)=-2^x x(\log 2)^2/(2^x-1+2^x x\log 2),$ and
%$$(u'(x)/v'(x))'=-(\log 2)^2\frac{2^x(2^x-x \log (2)-1)}{(2^x+2^x x \log (2)-1)^2}.$$
%Since, for $x>0$, $2^x=\sum_{n=0}^\infty (x\log 2)^n/n!>x \log (2)+1$. Thus, $u'(x)/v'(x)$ is strictly decreasing and similarly $u(x)/v(x)$. Therefore, $\varphi$ is strictly concave and $\log G_1(x)= \varphi(x)+\varphi(1-x)$ is strictly concave. Hence, $(\log G_1)'(x)= \varphi'(x)-\varphi'(1-x)$ is strictly decreasing
%and since, $(\log G_1)'(1/2)=0$. Then $\log G_1$ and $G_1$ are strictly increasing on $(0,1/2)$ and strictly decreasing on $(1/2,1)$.
\begin{corollary} \
\begin{enumerate}
\item For all $b>1$ $a\geq 0$ and $x>1$,  $$\frac{\zeta'(x)}{\zeta(x)}>\frac{b}{1-x}-\frac ax.$$
\item For all $a\geq 0$ and $x>1$
    $$\frac{\zeta'(x)}{\zeta(x)}<\frac{1}{1-x}-\frac ax+\gamma+a,$$

\item  For all $x\in(0,\infty)$
 $$\log(2\pi)-1+\frac1{1-x}>\frac{\zeta'(x)}{\zeta(x)}>\frac{1}{1-x}.$$
 The constant upper bound $\log(2\pi)-1$ and $a+\gamma$ are sharp.
 \end{enumerate}
\end{corollary}

  \begin{proposition}\label{pr3}  \

  \begin{enumerate}
  \item The function $u(s)=\zeta(s)\zeta(1/s)$ is strictly decreasing on $(0,1)$ and strictly increasing on $(1,+\infty)$ and for all $s>0$, $s\neq 1$
  $$\zeta(s)\zeta(1/s)<-\frac 12$$
  \item The function $h(s)=\zeta(1+s)\zeta(1-s)$ is strictly increasing on $(0,1)$ and
  $$\zeta(1+s)\zeta(1-s)<-\frac{\pi^2}{12},$$
  the upper bound is sharp.

  \end{enumerate}
  \end{proposition}

  {\bf Proof.}
a) Differentiate yields
   $$u'(s)=\zeta'(s)\zeta(1/s)-\frac1{s^2}\zeta(s)\zeta'(1/s),$$
   For $s>1$, we have $\zeta'(s)<0$ and from Remark \ref{r2} $\zeta(1/s)<0$ and $\zeta'(1/s)<0$. Therefore, $g'(s)>0$ and $g$ is strictly increasing on $(1,+\infty)$. From the functional equation $g(1/s)=g(s)$ we deduce that $g$ is strictly decreasing on $(0,1)$. Furthermore, by analytic continuation we have $\zeta(0)=-1/2$ and  $\lim_{s\to+\infty}\zeta(s)=1$, then
   $\lim_{s\to 0}g(s)=\zeta(0)=-1/2$.

  b) Differentiate yields
  $$h'(s)=\zeta'(1+s)\zeta(1-s)-\zeta(1+s)\zeta'(1-s).$$
  Since, for $s\in(0,1)$, $\zeta'(1+s)<0$ and by Remark \ref{r2}, $\zeta(1-s)<0$ and $\zeta'(1-s)<0$, then $h'(s)>0$ for all $s\in(0,1)$.

  Therefore, $h$ is strictly increasing on $(0,1)$. Furthermore, and by using $\zeta(0)=-1/2$ and $\zeta(2)=\pi^2/6$ we get the desired result.

\section{Proof of Theorem \ref{t1}}
From (a) of Proposition \ref{pr3} the expression of Theorem \ref{t1} is defined for all positive $x$.
The function $s\mapsto\zeta(s)\zeta(1/s)$ is strictly decreasing and negative on $(0,1)$ then $s\mapsto-1/(\zeta(s)\zeta(1/s))$ is strictly decreasing and positive on $(0,1)$, moreover, the function $s\mapsto\zeta(s)+\zeta(1/s)$ is strictly decreasing and positive on $(0,1)$ therefore, the function
$s\mapsto\zeta(s)+\zeta(1/s)/(\zeta(s)\zeta(1/s))$ is strictly increasing on $(0,1)$ and strictly decreasing on $(1,+\infty)$ by symmetry. Moreover,
$$\lim_{s\to 1}\frac{\zeta(s)+\zeta(1/s)}{\zeta(s)\zeta(1/s)}=0,$$
and
$$\lim_{s\to 0}\frac{\zeta(s)+\zeta(1/s)}{\zeta(s)\zeta(1/s)}=-1.$$
 This completes the proof of the theorem.

%\begin{acknowledgements}
%The author is very grateful to the referee for their several
%constructive comments of the first draft of this paper. %INCLUDE ANY ACKNOWLEDGEMENTS HERE
%\end{acknowledgements}
% IF YOU USE BIBTEX, THE FOLLOWING STYLE SHOULD BE USED
%\bibliographystyle{kjmath}
% INSERT THE NAME OF YOUR BIB FILE BELOW
%\bibliography{bibliography_file}

% BEFORE SUBMITTING YOUR MANUSCRIPT PLEASE COPY THE CONTENT OF THE BBL FILE THAT LATEX GENERATES HERE (IT WILL BE IN THE SAME DIRECTORY AS YOUR LATEX FILE) AND COMMENT THE PREVIOUS TWO LINES (\bibliographystyle and \bibliography), FOR EXAMPLE:

\end{document}